\documentclass{amsbook}
\include {mak}
\newcommand{\lbl}[1]{\label{#1}}
\newtheorem{theo}{Theorem}
\newtheorem{lemma}[theo]{Lemma}

\newtheorem{col} {Corollary}
\theoremstyle{definition}
\newtheorem{definition}[theo]{Definition}

\theoremstyle{remark}
\newtheorem{remark} {Remark}
\newcommand{\be}{\begin{equation}}
\newcommand{\ee}{\end{equation}}
\newcommand\bes{\begin{eqnarray}} \newcommand\ees{\end{eqnarray}}
\newcommand{\bess}{\begin{eqnarray*}}
\newcommand{\eess}{\end{eqnarray*}}
\numberwithin{equation}{section}

\begin{document}


\oddsidemargin 15.5mm \evensidemargin 15.5mm

\thispagestyle{plain}


\vspace{3cc}

\begin{center}
 {\Large\bf  New Generalization of Perturbed Ostrowski Type Inequalities and Applications
\rule{0mm}{6mm}\renewcommand{\thefootnote}{}
\footnotetext{

Keywords and Prases: perturbed Ostrowski type inequality, functions
of Lipschitzian type, numerical integration, cumulative distribution
functions. }}\vspace{1cc}

 {\large \textsc{Wen-jun Liu, Qiao-ling Xue, Jian-wei Dong } }\\

\end{center}

\vspace{1cc} \noindent \textsc{Abstract}: Generalizations of
Ostrowski type inequality for functions of Lipschitzian type are
established. Applications in numerical integration and  cumulative
distribution functions are also given.\vspace{3cc}

%
%

\section{Introduction}

In \cite{nu2004} N. Ujevi{\rm$\acute{c}$}  obtained the following
perturbed Ostrowski type inequality.

 \begin{theo}\lbl{t1}
Let $I\subset R $ be an open interval such that $[a,b] \subset I$
and let $f: I\rightarrow R $ be a  differentiable function such that
$\gamma\leq f'(t)\leq \Gamma, \ \forall t\in [a,b]$, for some
constants $\gamma, \Gamma\in R$. Then we have \be\begin{array}{lll}
\medskip
\ \ \
\displaystyle\left|(b-a)\left\{\left[f(x)-\frac{\Gamma+\gamma}{2}\left(x-\frac{a+b}{2}\right)\right](1-h)
+\frac{f(a)+f(b)}{2}h\right\}
-\int_a^bf(t)dt\right|\\
\leq \displaystyle
\frac{1}{2}\left[\frac{(b-a)^2}{4}[h^2+(h-1)^2]+\left(x-\frac{a+b}{2}\right)^2\right](\Gamma-\gamma),
\end{array}\label{1.1}\ee
where $a+h((b-a)/2)\leq x\leq b-h((b-a)/2)$ and $h\in [0,1].$
\end{theo}

In \cite{nu2002}, the same author proved the next result.

 \begin{theo}\lbl{t2}
Let the assumotions of Theorem \ref{t1} hold. Then for all
$a+h((b-a)/2)\leq x\leq b-h((b-a)/2)$ and $h\in [0,1]$, we have
\be\begin{array}{lll}
\medskip
\ \ \
\displaystyle\left|(b-a)\left\{\left[f(x)-\gamma\left(x-\frac{a+b}{2}\right)\right](1-h)
+\frac{f(a)+f(b)}{2}h\right\}
-\int_a^bf(t)dt\right|\\
\leq \displaystyle (b-a)\max \left\{h\frac{b-a}{2},
x-a-h\frac{b-a}{2}, b-x-h\frac{b-a}{2}\right\}(S-\gamma),
\end{array}\label{1.2}\ee
and \be\begin{array}{lll}
\medskip
\ \ \
\displaystyle\left|(b-a)\left\{\left[f(x)-\Gamma\left(x-\frac{a+b}{2}\right)\right](1-h)
+\frac{f(a)+f(b)}{2}h\right\}
-\int_a^bf(t)dt\right|\\
\leq \displaystyle (b-a)\max \left\{h\frac{b-a}{2},
x-a-h\frac{b-a}{2}, b-x-h\frac{b-a}{2}\right\}(\Gamma-S),
\end{array}\label{1.3}\ee
where $S=(f(b)-f(a))/(b-a).$
\end{theo}

All (\ref{1.1})-(\ref{1.3}) have been used to get the tighter error
bounds for the midpoint, trapezoid, and Simpson quadrature formulas
in numerical integration, respectively.

In this paper, we shall generalize Theorems \ref{t1} and \ref{t2} to
functions of some larger classes. Applications in numerical
integration and  cumulative distribution functions are also given.
For convenience, we define functions of Lipschitzian type as
follows:

 \begin{definition}\lbl{d3} The function $f: [a, b]\rightarrow R$ is said to be $L$-Lipschitzian on $[a,
 b]$
if for some $L>0$ and all $x, y\in [a, b]$,
$$|f(x)-f(y)|\leq L |x-y|.$$
 \end{definition}

  \begin{definition}\lbl{d4} The function $f: [a, b]\rightarrow R$ is said to be $(l, L)$-Lipschitzian on $[a,
 b]$
if
$$l(x_2-x_1)\leq f(x_2)-f(x_1)\leq L(x_2-x_1)\ \ \mbox{for}\ \ a\leq x_1\leq x_2\leq b,$$
where $l, L\in R$ with $l<L$.
 \end{definition}

 We will need the following well-known results.
 \begin{lemma}{( \cite{dac})} \lbl{l5} Let $h, g: [a, b]\rightarrow R$ be such that h is Riemann integrable
on $[a, b]$ and $g$ is $L$-Lipschitzian on $[a, b]$. Then
$$\left|\int_a^bh(t)dg(t)\right|\leq L\int_a^b|h(t)|dt.$$
  \end{lemma}
  \begin{lemma}{( \cite{dac})} \lbl{l6} Let $h, g: [a, b]\rightarrow R$ be such that h is continuous
on $[a, b]$ and $g$ is of bounded variation on $[a, b]$. Then
$$\left|\int_a^bh(t)dg(t)\right|\leq \max_{t\in [a,b]}|h(t)|\bigvee_a^b(g).$$
  \end{lemma}

%
\section{Main results}
\setcounter{equation}{0}
Our main results are as follows.

 \begin{theo}\lbl{t7}
Let $f: [a, b]\rightarrow R$   be $(l, L)$-Lipschitzian on $[a,
 b]$. Then we have \be\begin{array}{lll}
\medskip
\ \ \
\displaystyle\left|(b-a)\left\{\left[f(x)-\frac{L+l}{2}\left(x-\frac{a+b}{2}\right)\right](1-h)
+\frac{f(a)+f(b)}{2}h\right\}
-\int_a^bf(t)dt\right|\\
\leq \displaystyle
\frac{1}{2}\left[\frac{(b-a)^2}{4}[h^2+(h-1)^2]+\left(x-\frac{a+b}{2}\right)^2\right](L-l),
\end{array}\label{2.1}\ee
\be\begin{array}{lll}
\medskip
\ \ \
\displaystyle\left|(b-a)\left\{\left[f(x)-l\left(x-\frac{a+b}{2}\right)\right](1-h)
+\frac{f(a)+f(b)}{2}h\right\}
-\int_a^bf(t)dt\right|\\
\leq \displaystyle (b-a)\max \left\{h\frac{b-a}{2},
x-a-h\frac{b-a}{2}, b-x-h\frac{b-a}{2}\right\}(S-l),
\end{array}\label{2.2}\ee
and \be\begin{array}{lll}
\medskip
\ \ \
\displaystyle\left|(b-a)\left\{\left[f(x)-L\left(x-\frac{a+b}{2}\right)\right](1-h)
+\frac{f(a)+f(b)}{2}h\right\}
-\int_a^bf(t)dt\right|\\
\leq \displaystyle (b-a)\max \left\{h\frac{b-a}{2},
x-a-h\frac{b-a}{2}, b-x-h\frac{b-a}{2}\right\}(L-S),
\end{array}\label{2.3}\ee
for all $a+h((b-a)/2)\leq x\leq b-h((b-a)/2)$ and $h\in [0,1]$,
where $S=(f(b)-f(a))/(b-a).$
\end{theo}

\begin{proof}
Let $p: [a,b]^2 \rightarrow R$ be given by \be p(x,t):=\left\{
\begin{array}{lll}
\medskip t-[a+h\frac{b-a}{2}],\ &t\in [a,x]\\
\medskip t-[b-h\frac{b-a}{2}],\ &t\in (x,b],
\end{array} \right.\label{2.4}\ee
Put \be g(t):=f(t)-\frac{L+l}{2}t.\label{2.5}\ee It is easy to find
that the function $g: [a, b]\rightarrow R$ is $M$-Lipschitzian on
$[a, b]$ with $ M=\frac{L-l}{2}$. So, the Riemann-Stieltjes integral
$\int_a^b p(x, t) dg(t)$ exists. Using the integration by parts
formula for Riemann-Stieltjes integral, we have \be
\begin{array}{lll}
\medskip \displaystyle\int_a^bp(x,t)d g(t)\hskip -3pt
 &\displaystyle=\int_a^x\left(t-\left[a+h\frac{b-a}{2}\right]\right)dg(t)
+ \int_x^b\left(t-\left[b-h\frac{b-a}{2}\right]\right)dg(t)\\
\medskip &\displaystyle=(b-a)\left[g(x)(1-h)+\frac{g(a)+g(b)}{2}h\right]-\int_a^bg(t)dt.
\end{array} \label{2.6}\ee
By Lemma \ref{l5} we have \be
\left|(b-a)\left[g(x)(1-h)+\frac{g(a)+g(b)}{2}h\right]-\int_a^bg(t)dt\right|\leq
\frac{L-l}{2}\int_a^b|p(x,t)|dt.\label{2.7}\ee It is not difficult
to find that ( see \cite{dcr} ) \be
\int_a^b|p(x,t)|dt=\frac{(b-a)^2}{4}[h^2+(h-1)^2]+\left(x-\frac{a+b}{2}\right)^2,\label{2.8}\ee
and so from (\ref{2.7}) and (\ref{2.8}) we get \be\begin{array}{lll}
\medskip
\ \ \ \displaystyle\left|(b-a)\left[g(x)(1-h)
+\frac{g(a)+g(b)}{2}h\right]
-\int_a^bg(t)dt\right|\\
\leq \displaystyle
\frac{1}{2}\left[\frac{(b-a)^2}{4}[h^2+(h-1)^2]+\left(x-\frac{a+b}{2}\right)^2\right](L-l).
\end{array}\label{2.9}\ee
Consequently, the inequality (\ref{2.1}) follows from substituting
(\ref{2.5}) to the left hand side of the inequality (\ref{2.9}).

 Now
we proceed to prove the inequalities (\ref{2.2}) and (\ref{2.3}).
Put \be g_1(t):=f(t)-lt\ \ \mbox{and}\ \
g_2(t):=f(t)-Lt.\label{2.10}\ee It is easy to find that both $g_1,
g_2: [a, b]\rightarrow R$ are functions of bounded variation on
$[a,b]$ with \be \bigvee_a^b(g_1)=f(b)-f(a)-l(b-a)\ \ \mbox{and}\ \
\bigvee_a^b(g_2)=L(b-a)-[f(b)-f(a)].\label{2.11}\ee So, the
Riemann-Stieltjes integrals $\int_a^b p(x, t) dg_1(t)$  and
$\int_a^b p(x, t) dg_2(t)$  exist. Using the integration by parts
formula for Riemann-Stieltjes integral, we have \be \int_a^bp(x,t)d
g_1(t)=(b-a)\left[g_1(x)(1-h)+\frac{g_1(a)+g_1(b)}{2}h\right]-\int_a^bg_1(t)dt,
\label{2.12}\ee and
 \be \int_a^bp(x,t)d
g_2(t)=(b-a)\left[g_2(x)(1-h)+\frac{g_2(a)+g_2(b)}{2}h\right]-\int_a^bg_2(t)dt.
\label{2.13}\ee Then by  Lemma \ref{l6} we can deduce that $$
\left|(b-a)\left[g_1(x)(1-h)+\frac{g_1(a)+g_1(b)}{2}h\right]-\int_a^bg_1(t)dt\right|\leq
\max_{t\in [a,b]}|p(x,t)|\bigvee_a^b(g_1)$$ and
$$
\left|(b-a)\left[g_2(x)(1-h)+\frac{g_2(a)+g_2(b)}{2}h\right]-\int_a^bg_2(t)dt\right|\leq
\max_{t\in [a,b]}|p(x,t)|\bigvee_a^b(g_2).$$ Notice that
$$\max_{t\in [a,b]}|p(x,t)|=\max\left\{h\frac{b-a}{2},
x-a-h\frac{b-a}{2}, b-x-h\frac{b-a}{2}\right\},$$ and from
(\ref{2.11}), we get \be\begin{array}{lll}
\medskip
\ \ \ \displaystyle\left|(b-a)\left[g_1(x)(1-h)
+\frac{g_1(a)+g_1(b)}{2}h\right]
-\int_a^bg_1(t)dt\right|\\
\leq \displaystyle
 \displaystyle (b-a)\max \left\{h\frac{b-a}{2},
x-a-h\frac{b-a}{2}, b-x-h\frac{b-a}{2}\right\}(S-l),
\end{array}\label{2.14}\ee
and \be\begin{array}{lll}
\medskip
\ \ \ \displaystyle\left|(b-a)\left[g_2(x)(1-h)
+\frac{g_2(a)+g_2(b)}{2}h\right]
-\int_a^bg_2(t)dt\right|\\
\leq \displaystyle
 \displaystyle (b-a)\max \left\{h\frac{b-a}{2},
x-a-h\frac{b-a}{2}, b-x-h\frac{b-a}{2}\right\}(L-S),
\end{array}\label{2.15}\ee
where $S=(f(b)-f(a))/(b-a)$.

Consequently, inequalities (\ref{2.2}) and (\ref{2.3}) follow from
substituting (\ref{2.10}) to the left hand sides of (\ref{2.14}) and
(\ref{2.15}), respectively.
 \end{proof}

 \begin{col}\lbl{c1}
Under the assumptions of Theorem \ref{t7}, we have
\be\begin{array}{lll}
\medskip
\ \ \
\displaystyle\left|(b-a)\left[f(x)-\frac{L+l}{2}\left(x-\frac{a+b}{2}\right)\right]
-\int_a^bf(t)dt\right|\\
\leq \displaystyle
\frac{1}{2}\left[\frac{(b-a)^2}{4}+\left(x-\frac{a+b}{2}\right)^2\right](L-l),
\end{array}\label{2.16}\ee
\be\begin{array}{lll}
\medskip
\ \ \
\displaystyle\left|(b-a)\left[f(x)-l\left(x-\frac{a+b}{2}\right)\right]
-\int_a^bf(t)dt\right|\\
\leq \displaystyle
(b-a)\left[\frac{b-a}{2}+\left|x-\frac{a+b}{2}\right|\right](S-l),
\end{array}\label{2.17}\ee
and \be\begin{array}{lll}
\medskip
\ \ \
\displaystyle\left|(b-a)\left[f(x)-L\left(x-\frac{a+b}{2}\right)\right]
-\int_a^bf(t)dt\right|\\
\leq \displaystyle
(b-a)\left[\frac{b-a}{2}+\left|x-\frac{a+b}{2}\right|\right](L-S).
\end{array}\label{2.18}\ee
\end{col}
\begin{proof}
We set $h=0$ in the above theorem and utilize \be \max\{x-a,
b-x\}=\frac{1}{2}[b-a+|2x-a-b|]=\frac{b-a}{2}+\left|x-\frac{a+b}{2}\right|.\label{2.19}\ee
\end{proof}

 \begin{remark}\lbl{r1}
If we set $x=(a+b)/2$ in Corollary \ref{c1}, then we get
corresponding mid-point inequalities.
\end{remark}

 \begin{col}\lbl{c2}
Under the assumptions of Theorem \ref{t7}, we have \be
\left|\frac{b-a}{2}[f(a)+f(b)]-\int_a^bf(t)dt\right| \leq
\frac{(b-a)^2}{8}(L-l), \label{2.20}\ee \be
\left|\frac{b-a}{2}[f(a)+f(b)]-\int_a^bf(t)dt\right| \leq
\frac{(b-a)^2}{2}(S-l), \label{2.21}\ee and \be
\left|\frac{b-a}{2}[f(a)+f(b)]-\int_a^bf(t)dt\right| \leq
\frac{(b-a)^2}{2}(L-S), \label{2.22}\ee
\end{col}
\begin{proof}
We set $h=1$ in the above theorem and utilize \be
\max\left\{\frac{b-a}{2}, x-\frac{a+b}{2},
\frac{a+b}{2}-x\right\}=\frac{b-a}{2}.\label{2.23}\ee
\end{proof}

 \begin{col}\lbl{c3}
Under the assumptions of Theorem \ref{t7}, we have
\be\begin{array}{lll}
\medskip
\ \ \
\displaystyle\left|(b-a)\left[\frac{1}{2}f(x)-\frac{L+l}{4}\left(x-\frac{a+b}{2}\right)+\frac{f(a)+f(b)}{4}\right]
-\int_a^bf(t)dt\right|\\
\leq \displaystyle
\frac{1}{2}\left[\frac{(b-a)^2}{8}+\left(x-\frac{a+b}{2}\right)^2\right](L-l),
\end{array}\label{2.24}\ee
\be\begin{array}{lll}
\medskip
\ \ \
\displaystyle\left|(b-a)\left[\frac{1}{2}f(x)-\frac{l}{2}\left(x-\frac{a+b}{2}\right)+\frac{f(a)+f(b)}{4}\right]
-\int_a^bf(t)dt\right|\\
\leq \displaystyle
(b-a)\left[\frac{b-a}{4}+\left|x-\frac{a+b}{2}\right|\right](S-l),
\end{array}\label{2.25}\ee
and \be\begin{array}{lll}
\medskip
\ \ \
\displaystyle\left|(b-a)\left[\frac{1}{2}f(x)-\frac{L}{2}\left(x-\frac{a+b}{2}\right)+\frac{f(a)+f(b)}{4}\right]
-\int_a^bf(t)dt\right|\\
\leq \displaystyle
(b-a)\left[\frac{b-a}{4}+\left|x-\frac{a+b}{2}\right|\right](L-S).
\end{array}\label{2.26}\ee
\end{col}
\begin{proof}
We set $h=1/2$ in the above theorem and utilize \be
\max\left\{\frac{b-a}{4}, x-\frac{3a+b}{4},  \frac{a+3b}{4}-x
\right\}=\frac{b-a}{4}+\left|x-\frac{a+b}{2}\right|.\label{2.27}\ee
\end{proof}

 \begin{remark}\lbl{r2}
If we set $x=(a+b)/2$ in Corollary \ref{c3}, then we get
corresponding three point inequalities ( i.e. the average of a
mid-point and trapezoid type rules ).
\end{remark}

 \begin{col}\lbl{c4}
Under the assumptions of Theorem \ref{t7}, we have
\be\begin{array}{lll}
\medskip
\ \ \
\displaystyle\left|\frac{b-a}{6}[f(a)+4f(x)+f(b)]-\frac{L+l}{3}\left(x-\frac{a+b}{2}\right)
-\int_a^bf(t)dt\right|\\
\leq \displaystyle
\frac{1}{2}\left[\frac{5}{36}(b-a)^2+\left(x-\frac{a+b}{2}\right)^2\right](L-l),
\end{array}\label{2.28}\ee
\be\begin{array}{lll}
\medskip
\ \ \
\displaystyle\left|\frac{b-a}{6}[f(a)+4f(x)+f(b)]-\frac{2l}{3}\left(x-\frac{a+b}{2}\right)
-\int_a^bf(t)dt\right|\\
\leq \displaystyle
(b-a)\left[\frac{b-a}{3}+\left|x-\frac{a+b}{2}\right|\right](S-l),
\end{array}\label{2.29}\ee
and \be\begin{array}{lll}
\medskip
\ \ \
\displaystyle\left|\frac{b-a}{6}[f(a)+4f(x)+f(b)]-\frac{2L}{3}\left(x-\frac{a+b}{2}\right)
-\int_a^bf(t)dt\right|\\
\leq \displaystyle
(b-a)\left[\frac{b-a}{3}+\left|x-\frac{a+b}{2}\right|\right](L-S).
\end{array}\label{2.30}\ee
\end{col}
\begin{proof}
We set $h=1/3$ in the above theorem and utilize \be
\max\left\{\frac{b-a}{6}, x-\frac{5a+b}{6},
\frac{a+5b}{6}-x\right\}=\frac{b-a}{3}+\left|x-\frac{a+b}{2}\right|.\label{2.31}\ee
\end{proof}

 \begin{remark}\lbl{r3}
If we set $x=(a+b)/2$ in Corollary \ref{c4}, then we get
corresponding Simpson's inequalities.\end{remark}

 \begin{remark}\lbl{r4}
It is interesting to note that the smallest bound for
(\ref{2.1})-(\ref{2.3}) is obtained at $h=1/2$ for fixed $x$. Thus
the quadrature rules (\ref{2.24})-(\ref{2.26}) comprised of the
linear combination of the perturbed mid-point and trapezoidal rules
are optimal and has a lower bound than the perturbed Simpson's rules
(\ref{2.28})-(\ref{2.30}).\end{remark}

\begin{remark}\lbl{r5}
It is clear that Theorem \ref{t7} can be regarded as generalization
of Theorems \ref{t1} and  \ref{t2}.
\end{remark}

\begin{theo}\lbl{t8}
Let $f: [a, b]\rightarrow R$   be $L$-Lipschitzian on $[a,
 b]$. Then we have \be\begin{array}{lll}
\medskip
\ \ \ \displaystyle\left|(b-a)\left[f(x)(1-h)
+\frac{f(a)+f(b)}{2}h\right]
-\int_a^bf(t)dt\right|\\
\leq \displaystyle
\left[\frac{(b-a)^2}{4}[h^2+(h-1)^2]+\left(x-\frac{a+b}{2}\right)^2\right]L,
\end{array}\label{2.32}\ee
and \be\begin{array}{lll}
\medskip
\ \ \
\displaystyle\left|(b-a)\left\{\left[f(x)+L\left(x-\frac{a+b}{2}\right)\right](1-h)
+\frac{f(a)+f(b)}{2}h\right\}
-\int_a^bf(t)dt\right|\\
\leq \displaystyle (b-a)\max \left\{h\frac{b-a}{2},
x-a-h\frac{b-a}{2}, b-x-h\frac{b-a}{2}\right\}(S+L),
\end{array}\label{2.33}\ee
for all $a+h((b-a)/2)\leq x\leq b-h((b-a)/2)$ and $h\in [0,1]$,
where $S=(f(b)-f(a))/(b-a).$
\end{theo}
\begin{proof}
We get inequality (\ref{2.32}) and (\ref{2.33}) immediately by
taking $l=-L$ in (\ref{2.1}) and (\ref{2.2}).
\end{proof}

\section{Applications in numerical integration} \setcounter{equation}{0}
We restrict further considerations to the perturbed three point
rules. We also emphasize that similar considerations can be done for
all quadrature rules considered in the previous section.

\begin{theo}\lbl{t9}
Let all the assumptions of Theorem \ref{t7} hold. If
$I_n=\{a=x_0<x_1<\cdots<x_n=b\}$ is a given subdivision of the
interval $[a, b]$ and $h_i=x_{i+1}-x_i, i=0,1,2,\cdots, n-1$, then

\be \label{3.1}\int_a^bf(t)dt=A_{Ll}(I_n,\xi,f)+R_{Ll}(I_n,\xi,f),
\ee where \be
\label{3.2}A_{Ll}(I_n,\xi,f)=\frac{1}{2}\sum_{i=0}^{n-1}f(\xi_i)h_i
+\frac{1}{2}\sum_{i=0}^{n-1}\frac{f(x_i)+f(x_{i+1})}{2}h_i
-\frac{L+l}{4}
\sum_{i=0}^{n-1}\left(\xi_i-\frac{x_i+x_{i+1}}{2}\right)h_i, \ee for
$x_i\leq \xi_i\leq x_{i+1}, i=0, 1, 2,\cdots, n-1.$ The remainder
term satisfies\be \label{3.3}|R_{Ll}(I_n,\xi,f)|\leq
\sum_{i=0}^{n-1}\frac{1}{2}\left[\frac{h_i^2}{8}+\left(\xi_i-\frac{x_i+x_{i+1}}{2}\right)^2\right](L-l).
\ee

Also, \be \label{3.4}\int_a^bf(t)dt=A_l(I_n,\xi,f)+R_{l}(I_n,\xi,f),
\ee where \be
\label{3.5}A_l(I_n,\xi,f)=\frac{1}{2}\sum_{i=0}^{n-1}f(\xi_i)h_i
+\frac{1}{2}\sum_{i=0}^{n-1}\frac{f(x_i)+f(x_{i+1})}{2}h_i
-\frac{l}{2}
\sum_{i=0}^{n-1}\left(\xi_i-\frac{x_i+x_{i+1}}{2}\right)h_i, \ee
and\be \label{3.6}|R_{l}(I_n,\xi,f)|\leq
\sum_{i=0}^{n-1}h_i\left[\frac{h_i}{4}+\left|\xi_i-\frac{x_i+x_{i+1}}{2}\right|\right](S_i-l),
\ee where $S_i=f(x_{i+1}-f(x_i))/h_i, i=0, 1, 2,\cdots, n-1.$

Also, \be \label{3.7}\int_a^bf(t)dt=A_L(I_n,\xi,f)+R_{L}(I_n,\xi,f),
\ee where \be
\label{3.8}A_L(I_n,\xi,f)=\frac{1}{2}\sum_{i=0}^{n-1}f(\xi_i)h_i
+\frac{1}{2}\sum_{i=0}^{n-1}\frac{f(x_i)+f(x_{i+1})}{2}h_i
-\frac{L}{2}
\sum_{i=0}^{n-1}\left(\xi_i-\frac{x_i+x_{i+1}}{2}\right)h_i, \ee and
\be \label{3.9}|R_{L}(I_n,\xi,f)|\leq
\sum_{i=0}^{n-1}h_i\left[\frac{h_i}{4}+\left|\xi_i-\frac{x_i+x_{i+1}}{2}\right|\right](L-S_i).
\ee
\end{theo}
\begin{proof}
Apply Corollary \ref{c3} to the interval $[x_i, x_{i+1}],
i=0,1,2,\cdots, n-1 $ and sum. Then use the triangle inequality to
obtain the desired result.
 \end{proof}
\begin{remark}\lbl{r6}
If we set $\xi_i=(x_i+x_{i+1})/2$ in Theorem \ref{t9}, then we get
corresponding composite rules which do not depend on $\xi_i$.
\end{remark}
\begin{remark}\lbl{r7}
Note that we can apply quadrature rules in \cite{nu2004} and
\cite{nu2002} only if $f\in C^1[a, b]$, while we can apply here if
$f$ is $(l, L)$-Lipschitzian. Hence, the above obtained result
enlarges the applicability of the quadrature rules.
\end{remark}

\section{Applications  for cumulative distribution functions} \setcounter{equation}{0}

Now we consider some applications for cumulative distribution
functions.

 Let $X$ be a random variable having the probability density
function $f: [a, b]\rightarrow R_+$ and the cumulative distribution
function $F(x)=Pr(X\leq x),$ i.e.,
$$F(x)=\int_a^xf(t)dt, \ \ x\in [a,b].$$
$E(X)$ is the expectation of $X$. Then we have the following
inequality.
\begin{theo}\lbl{t10}
With the above assumptions and if there exist constants $M,m$ such
that $0\leq m\leq f(t)\leq M$ for all $t\in [a, b]$, then we have
the inequalities \be\begin{array}{lll}
\medskip
\ \ \ \displaystyle\left|(b-a)\left\{\left[P_r(X\leq
x)-\frac{M+m}{2}\left(x-\frac{a+b}{2}\right)\right](1-h)
+\frac{1}{2}h\right\}
-(b-E(X))\right|\\
\leq \displaystyle
\frac{1}{2}\left[\frac{(b-a)^2}{4}[h^2+(h-1)^2]+\left(x-\frac{a+b}{2}\right)^2\right](M-m),
\end{array}\label{4.1}\ee
\be\begin{array}{lll}
\medskip
\ \ \ \displaystyle\left|(b-a)\left\{\left[P_r(X\leq
x)-m\left(x-\frac{a+b}{2}\right)\right](1-h) +\frac{1}{2}h\right\}
-(b-E(X))\right|\\
\leq \displaystyle (b-a)\max \left\{h\frac{b-a}{2},
x-a-h\frac{b-a}{2},
b-x-h\frac{b-a}{2}\right\}\left(\frac{1}{b-a}-m\right),
\end{array}\label{4.2}\ee
and \be\begin{array}{lll}
\medskip
\ \ \ \displaystyle\left|(b-a)\left\{\left[P_r(X\leq
x)-M\left(x-\frac{a+b}{2}\right)\right](1-h) +\frac{1}{2}h\right\}
-(b-E(X))\right|\\
\leq \displaystyle (b-a)\max \left\{h\frac{b-a}{2},
x-a-h\frac{b-a}{2},
b-x-h\frac{b-a}{2}\right\}\left(M-\frac{1}{b-a}\right),
\end{array}\label{4.3}\ee
for all $a+h((b-a)/2)\leq x\leq b-h((b-a)/2)$ and $h\in [0,1]$.
\end{theo}
\begin{proof}
It is easy to show that the function $F(x)=\int_a^xf(t)dt$ is $(m,
M)$-Lipschitzian on $[a, b]$. So, by Theorem \ref{t7} we get
\be\begin{array}{lll}
\medskip
\ \ \
\displaystyle\left|(b-a)\left\{\left[F(x)-\frac{M+m}{2}\left(x-\frac{a+b}{2}\right)\right](1-h)
+\frac{F(a)+F(b)}{2}h\right\}
-\int_a^bF(t)dt\right|\\
\leq \displaystyle
\frac{1}{2}\left[\frac{(b-a)^2}{4}[h^2+(h-1)^2]+\left(x-\frac{a+b}{2}\right)^2\right](M-m),
\end{array}\label{4.4}\ee
\be\begin{array}{lll}
\medskip
\ \ \
\displaystyle\left|(b-a)\left\{\left[F(x)-m\left(x-\frac{a+b}{2}\right)\right](1-h)
+\frac{F(a)+F(b)}{2}h\right\}
-\int_a^bF(t)dt\right|\\
\leq \displaystyle (b-a)\max \left\{h\frac{b-a}{2},
x-a-h\frac{b-a}{2}, b-x-h\frac{b-a}{2}\right\}(S-m),
\end{array}\label{4.5}\ee
and \be\begin{array}{lll}
\medskip
\ \ \
\displaystyle\left|(b-a)\left\{\left[F(x)-M\left(x-\frac{a+b}{2}\right)\right](1-h)
+\frac{F(a)+F(b)}{2}h\right\}
-\int_a^bF(t)dt\right|\\
\leq \displaystyle (b-a)\max \left\{h\frac{b-a}{2},
x-a-h\frac{b-a}{2}, b-x-h\frac{b-a}{2}\right\}(M-S),
\end{array}\label{4.6}\ee where $S=(F(b)-F(a))/(b-a).$

As $F(a)=0, F(b)=1$, and
$$\int_a^bF(t)dt=b-E(X),$$
then we can easily deduce inequalities (\ref{4.1})-(\ref{4.3}).
 \end{proof}

  \begin{col}\lbl{c5}
Under the assumptions of Theorem \ref{t10}, we have
\be\begin{array}{lll}
\medskip
\ \ \ \displaystyle\left|(b-a)\left[\frac{1}{2}P_r(X\leq
x)-\frac{M+m}{4}\left(x-\frac{a+b}{2}\right)+\frac{1}{4}\right]
-(b-E(X))\right|\\
\leq \displaystyle
\frac{1}{2}\left[\frac{(b-a)^2}{8}+\left(x-\frac{a+b}{2}\right)^2\right](M-m),
\end{array}\label{4.7}\ee
\be\begin{array}{lll}
\medskip
\ \ \ \displaystyle\left|(b-a)\left[\frac{1}{2}P_r(X\leq
x)-\frac{m}{2}\left(x-\frac{a+b}{2}\right)+\frac{1}{4}\right]
-(b-E(X))\right|\\
\leq \displaystyle
(b-a)\left[\frac{b-a}{4}+\left|x-\frac{a+b}{2}\right|\right]\left(\frac{1}{b-a}-m\right),
\end{array}\label{4.8}\ee
and \be\begin{array}{lll}
\medskip
\ \ \ \displaystyle\left|(b-a)\left[\frac{1}{2}P_r(X\leq
x)-\frac{M}{2}\left(x-\frac{a+b}{2}\right)+\frac{1}{4}\right]
-(b-E(X))\right|\\
\leq \displaystyle
(b-a)\left[\frac{b-a}{4}+\left|x-\frac{a+b}{2}\right|\right]\left(M-\frac{1}{b-a}\right).
\end{array}\label{4.9}\ee
\end{col}

 \begin{proof}
We set $h=1/2$ in the above theorem.
\end{proof}

 \begin{col}\lbl{c6}
Under the assumptions of Theorem \ref{t10}, we have
\be\left|\left[\frac{1}{2}P_r\left(X\leq
\frac{a+b}{2}\right)+\frac{1}{4}\right] -\frac{b-E(X)}{b-a}\right|
\leq
 \frac{b-a}{16}(M-m),
\label{4.10}\ee \be\left|\left[\frac{1}{2}P_r\left(X\leq
\frac{a+b}{2}\right)+\frac{1}{4}\right] -\frac{b-E(X)}{b-a}\right|
\leq \frac{b-a}{4}\left(\frac{1}{b-a}-m\right), \label{4.11}\ee and
\be\left|\left[\frac{1}{2}P_r\left(X\leq
\frac{a+b}{2}\right)+\frac{1}{4}\right] -\frac{b-E(X)}{b-a}\right|
\leq
 \frac{b-a}{4}\left(M-\frac{1}{b-a}\right).
\label{4.12}\ee
\end{col}

 \begin{proof}
We set $x=\frac{a+b}{2}$ in Corollary \ref{c5}.
\end{proof}

\begin{col}\lbl{c7}
Under the assumptions of Theorem \ref{t10}, we have \be
\left|E(x)-\frac{a+3b}{4}+\frac{1}{8}(M+m)(b-a)^2\right|\leq
\frac{3}{16}(b-a)^2(M-m),\label{4.13}\ee
\be\left|E(x)-\frac{a+3b}{4}+\frac{1}{4}m(b-a)^2\right|\leq
\frac{3}{4}(b-a)^2(M-m),\label{4.14}\ee and
\be\left|E(x)-\frac{a+3b}{4}+\frac{1}{4}M(b-a)^2\right|\leq
\frac{3}{4}(b-a)^2(M-m).\label{4.15}\ee
\end{col}

 \begin{proof}
We set $x=a$ or $x=b$ in Corollary \ref{c5}.
\end{proof}

\begin{remark}\lbl{r8}
Similar results can be obtained when set $h=0, 1$ or $1/3$ in the
above theorem.
\end{remark}

{\bf Acknowledgements:} The first author was supported by the
Science Research Foundation of NUIST, and the third author was
supported by Youth Natural Science Foundation of Zhengzhou Institute
of Aeronautical Industry Management under Grant No.Q05K066.
\vspace{0.8cc}

\begin{center}{\small\bf REFERENCES}
\end{center}

\vspace{0.8cc}
\newcounter{ref}

\begin{list}{\small \arabic{ref}.}{\usecounter{ref} \leftmargin 4mm \itemsep  1mm}


\bibitem{bd}
N. S. Barnett And S. S. Dragomir, \newblock{\em Some inequalities
for probability, expectation and variance of random variable defined
over a finite interval}, { Computers and Mathematics with
Applications}, {\bf 43}(2002), 1319-1357.

\bibitem{cd}
P. Cerone and S. S. Dragomir, Trapezoidal-type Rules from an
Inequalities Point of View, Handbook of Analytic-Computational
Methods in Applied Mathematics, Editor: G. Anastassiou, CRC Press,
New York, (2000), 65-134.

\bibitem{cd1}
P. Cerone and S. S. Dragomir, Midpoint-type Rules from an
Inequalities Point of View, Handbook of Analytic-Computational
Methods in Applied Mathematics, Editor: G. Anastassiou, CRC Press,
New York, (2000), 65-134.

\bibitem{dac}
S. S. Dragomir, R. P. Agarwal and P. Cerone, \newblock{\em On
Simpson¡¯s inequality and applications},  { J. Inequal. Appl.} {\bf
5} (2000), 533-579.

\bibitem{dcr}
S. S. Dragomir, P. Cerone and J. Roumeliotis, \newblock{\em A new
generalization of Ostrowski's integral inequality for mappings whose
derivatives are bounded and applications in numerical integration
and for special means},  {  Appl. Math. Lett.} {\bf 13}(1) (2000),
19-25.

\bibitem{l}
W. J. Liu, C. C. Li and J. W. Dong, \newblock{\em Note on Qi's
Inequality and Bougoffa's Inequality},  {J. Inequal. Pure Appl.
Math.} {\bf 7}(4) (2006), Art.129.

\bibitem{l2}
W. J. Liu, Q. L. Xue and J. W. Dong, \newblock{\em New
Generalization of Perturbed Trapezoid and Mid-point Inequalities and
Applications}, (submitted).

\bibitem{l3}
W. J. Liu, J. Luo and J. W. Dong, \newblock{\em Further
Generalization of Perturbed Mid-point and Trapezoid Inequalities and
Applications}, (submitted).

\bibitem{lz}
Z. Liu, \newblock{\em Some Ostrowski-Gr$\ddot{u}$ss type
inequalities and applications}, { Computers and Mathematics with
Applications}, (2007), doi:10.1016/j.camwa.2006.12.021.

\bibitem{nu2004}
N. Ujevi{\rm$\acute{c}$},
\newblock{\em A Generalization of Ostrowski's Inequality
and Applications in Numerical Integration}, {Applied Mathematics
Letters 17}, (2004), 133-137.

\bibitem{nu2002} N.
Ujevi{\rm$\acute{c}$},
\newblock{\em Pertubations of an Ostrowski Type Inequality
and Applications},
 {IJMMS}, {\bf 32}(8) (2002),
491-500.

\bibitem{nu5}
N. Ujevi{\rm$\acute{c}$}, \newblock{\em Error inequalities for a
generalized trapezoid rule}, {Appl. Math. Lett.}, {\bf 19} (2006),
32-37.

\end{list}

\vspace{0.8cc}

College of Mathematics and Physics, Nanjing University of
Information Science and Technology, Nanjing 210044, China

{\small \noindent

 E--mail: wjliu@nuist.edu.cn}

\vspace{1cc} College of Mathematics and Physics, Nanjing University
of Information Science and Technology, Nanjing 210044, China

{\small \noindent

 E--mail: qlx\_\,1@yahoo.com.cn}

 \vspace{1cc} Department of Mathematics and Physics,
 Zhengzhou Institute of Aeronautical Industry Management,
Zhengzhou 450015, China.

 {\small \noindent

 E--mail: dongjianweiccm@163.com}

 \end{document}